%% file: doc.tex
\documentclass[11pt]{article}

\usepackage[a4paper, top=50pt, bottom=50pt]{geometry}
\usepackage[utf8]{inputenc}
\usepackage{amssymb, amsmath, xcolor, mathtools, makecell, tikz,
            multicol, microtype, xspace, mathabx, float, %mathspec, 
            hyperref}
\usepackage{enumitem}
\usepackage[lofdepth,lotdepth]{subfig}
\usetikzlibrary{automata}

\newcommand{\acknowledgements}{\subsection*{Acknowledgements}}
\newcommand*\getlength[1]{\number#1}
\newlength{\parindentS}
\usepackage[thmmarks,amsmath]{ntheorem}
\theoremstyle{plain}
\theoremseparator{:}
\newtheorem{defi}{Definition}
\newtheorem{thm}{Theorem}
\newtheorem{lemme}{Lemma}
\newtheorem{coro}{Corollary}
\newtheorem{prop}{Proposition}

\theorembodyfont{\upshape}
\newtheorem{ex}{Example}
\newtheorem{remarq}{Remark}

\theoremheaderfont{\sc}
\theoremstyle{nonumberplain} 
\theoremsymbol{\rule{1ex}{1ex}} 
\newtheorem{preuve}{Proof}
\newcommand{\NN}{\ensuremath{\mathbb{N}}}
\newcommand{\seq}{\ensuremath{\geqslant}}
\newcommand{\ieq}{\ensuremath{\leqslant}}
\newcommand{\card}[1]{\ensuremath{\# #1}}
\makeatletter
\newcommand*\oeis@star[1]{\textcolor{blue}{\href{http://oeis.org/#1}{\textsf{#1}}}}
\newcommand*\oeis@nostar[1]{[\oeis@star{#1}]}
\newcommand\oeis{\@ifstar{\oeis@star}{\oeis@nostar}}
\makeatother

\newcommand{\bij}{\ensuremath{\zeta}}

\newcommand{\seqPF}{\ensuremath{\chi}} 
\newcommand{\seqPFn}[1]{\ensuremath{\seqPF\left(#1\right)}} 
\newcommand{\PF}{\ensuremath{\mathcal{F}}}

\newcommand{\pf}{\ensuremath{\mathfrak{f}}}
\newcommand{\spf}{\ensuremath{\mathfrak{s}}}

\newcommand{\speciesSet}{\ensuremath{\mathtt{E}}}
\newcommand{\speciesOne}{\ensuremath{\boldsymbol{1}}}
\newcommand{\speciesA}{\ensuremath{\mathtt{P}}}
\newcommand{\speciesB}{\ensuremath{\mathtt{Q}}}

\newcommand{\ADFA}{\ensuremath{\mathcal{A}}}  %les ADFA
\newcommand{\niADFA}{\ensuremath{\mathcal{N}}} % non-initial ADFA : (F, \delta)
\newcommand{\niscADFA}{\ensuremath{\mathcal{S}}} % non-initial simple coreach ADFA
\renewcommand{\split}{\kappa} 
\newcommand{\eADFA}{\ensuremath{\mathcal{E}}} % (F, \delta) étendue
\newcommand{\mADFA}{\ensuremath{\mathcal{M}}} % (F, \delta) étendue
\newcommand{\esADFA}{\ensuremath{\widebar{\eADFA}}} % (F, \delta) étendue simple
\newcommand{\adfa}{\ensuremath{\mathfrak{a}}}
\newcommand{\dadfa}{\ensuremath{\mathfrak{d}}}

\newcommand{\eadfa}{\ensuremath{\mathfrak{e}}}

\newcommand{\madfa}{\ensuremath{\mathfrak{m}}}
\newcommand{\aut}{\ensuremath{\Theta}}
\newcommand{\autb}{\ensuremath{\Pi}}

\let\SOld\S
\renewcommand{\S}{\ensuremath{\mathbf{S}}}
\newcommand{\ncch}{\mathbf{ch}}

\newcommand*\pvsectiondefrec{\SOld 2.1} %section def récursive
\newcommand*\pveqrecursiveconstruction{2} %% définition récursive PF
\newcommand*\pvcorsolutionrecurrence{Proposition 2.5} %% solution de la
\newcommand*\pvthmFrobeniuscharacteristic{Theorem 3.4} %% le theoreme de la
\newcommand*\pvexplicationSchi{Eq. (7) and (9)}

\usepackage[affil-it]{authblk}

\author{Jean-Baptiste
Priez\thanks{{\url{jean-baptiste.priez@lri.fr}}}}
\affil{\href{http://www.lri.fr}{Laboratoire de Recherche en Informatique}\\
Université Paris-Sud}

\title{Enumeration of minimal acyclic automata
\textit{via} generalized parking functions}

\date{\today}

\definecolor{background}{RGB}{255, 255, 255}
\definecolor{fontcolor}{rgb}{0,0,0}

\everymath\expandafter{\the\everymath \color{gray!70!blue}}
\hypersetup{
 colorlinks=false,
 citecolor=Violet,
 linkcolor=Red,
 urlcolor=Blue
}

\begin{document}

\maketitle

\begin{abstract}\noindent
  \begin{description}
  \item[English] We give an exact enumerative formula for the minimal acyclic
    deterministic finite automata. This formula is obtained from a bijection
    between a family of generalized parking functions and the transitions
    functions of acyclic automata.
  \item[French] On donne une formule d'énumération exacte des automates finites
    déterministes acycliques minimaux. Cette formule s'obtient à partir d'une
    bijection entre une famille fonctions de parking généralisées et les
    fonctions de transitions des automates acycliques.
  \end{description}
\end{abstract}
\tableofcontents
\section*{Introduction} 
    \label{sec:intro} 
    \input{intro}
\section{Acyclic deterministic finite automata} 
    \label{sec:adfa}
    \input{sec_adfa}
\section{Generalized parking functions} 
    \label{sec:pf}
    \input{sec_pf} 
\section{Main result}
    \label{sec:main}
    \input{sec_main}    
%    \input{parking_functions}
%\input{pf}
%
%\section{Acyclic deterministic finite automata} \label{sec:adfa} \input{adfa}
%
\acknowledgements
\label{sec:ack}
I would like to thank anonymous reviewers for their precious comments
and also François \textsc{Bergeron} for the truly rewarding summer internship
in Montreal (supported by LIA LIRCO).
%
%\nocite{*}
%\bibliographystyle{abbrvnat}
% use the following instead if you encounter problems
\small 
\setlength{\parindent}{-20pt}
\bibliographystyle{alpha}
\bibliography{biblio}
\label{sec:biblio}

\end{document}

%% file: intro.tex
The study of the enumeration of minimal acyclic deterministic finite automata
(MADFA) has been undertaken by several authors in the last decade.
\textsc{Domaratzki} and \textit{al.} \cite{domaratzki2002number,
domaratzki2003improved, domaratzki2004combinatorial} presented different lower
and upper bounds, \textsc{Câmpeanu} and \textsc{Ho} \cite{campeanu2004maximum}
gave a good upper bound of MADFA with some constraints, \textsc{Almeida} and
\textit{al}. \cite{reis2005representation, almeida2007enumeration,
almeida2008exact} gave a canonical representation of MADFA and obtained a
method for exact generation.

In this paper, we refer to the study of \textsc{Liskovets}
\cite{liskovets2006exact}. The latter gave a recurrence relation to
enumerate acyclic finite deterministic automata (ADFA). The main idea of
\cite{liskovets2006exact} has been to define an extended notion of ADFA with
more than one \emph{absorbing state}. Unfortunately its approach of the
enumeration of ADFA is not fine enough to enumerate MADFA. The goal of this
paper is to give a finer enumeration of ADFA. In particular, the formula given
in this paper expresses properties on the right language of ADFA. Therefore, it
allows to enumerate MADFA. The main tool is a bijection with generalized
parking functions.

\textsc{Virmaux} and the author studied in \cite{priezvirmaux} the
generalization of parking functions defined by \cite{stanley2002polytope} and
gave a generalized generating series in non-commutative symmetric functions as
in \cite{novelli2008noncommutative}. Those generating series, called the
non-commutative Frobenius characteristic of the natural action of the $0$-Hecke
algebra, contain substantial information on combinatorial objects. 

In the first section, we recall the definition of extended ADFA of
\cite{liskovets2006exact} and then enrich the definition with some constraints.

In the second section, we recall the background on generalized parking functions
and we give an isomorphism with a noteworthy family of generalized parking
functions. We explicit the bijection and we transport some interesting
properties on the right language of ADFA directly on parking functions.
The substantial information provided by the Frobenius characteristic shows us
how to extract sub-families of parking functions. We use this information to
extract the sub-family of parking functions which exactly encodes (by the
bijection) the ADFA such that all states have their right language distinct.
Finally in the same way as \cite{liskovets2006exact}, one uses a
bijection between ADFA and couple of extended ADFA (with constraints) and MADFA.   
This defines a recurrence relation for the enumeration of minimal acyclic
deterministic finite automata.

%% file: sec_adfa.tex
For a basic background on automaton, the reader may consult
\cite{hopcroft1979introduction}. 

A \emph{deterministic finite automaton} (DFA) of $n$ states (labeled by a set
$N$) over an alphabet $\Sigma$ of $k$ symbols is a tuple $(i,A,\delta)$ with
\begin{itemize}
\begin{minipage}{0.4\linewidth}\color{fontcolor}
  \item $i \in N$ the \textit{initial state},
  \item $A \subset N$ the \textit{accepting states}, and
\end{minipage}
\begin{minipage}{0.5\linewidth}\color{fontcolor}
    \item $\delta : N \times \Sigma \to N \cup \{\emptyset\}$ the
  \textit{transition function}.\\[9pt]
\end{minipage}
\end{itemize}
% \begin{multicols}{2}
% \begin{itemize}
% \end{itemize}
% \end{multicols}
The special state $\emptyset$ is called the \emph{absorbing state}. We
consider $N$ the set of states to avoid worrying about ``well-labeled'' states
from $1$ to $n$.
%%
%\begin{multicols}{2}
\setlength{\parskip}{-3pt}
\paragraph{Extended transition function}  
We extend the transition function $\delta$  recursively to words
on $\Sigma$:
\begin{align*}
    \delta^*\; :\; N\cup\{\emptyset\} \times \Sigma^*
    \quad\longrightarrow\quad N \cup \{\emptyset\}\,,
\end{align*} 
by setting $\delta^*(q, aw) := \delta^*(\delta(q, a), w)$, for any $w \in
\Sigma^*$ and any $a \in \Sigma$; $\delta^*(q, \epsilon) = q$ (with $\epsilon$ the
empty word), for any state $q \in N$; and, $\delta^*(\emptyset, w) :=
\emptyset$, for any $w$.
\paragraph{Transitions of a state} We denote $\delta_q$ the
underlying transition function at $q$ defined by $\delta_q(a) := \delta(q, a)$.
\paragraph{Accepting status} The \emph{accepting status} of a state $q$ denotes
if $q$ is accepting or not (\textit{true} or \textit{false}).
\paragraph{Right language} The \emph{right language} of a state $q$ is the
language: $RL(q) := \{w \in \Sigma^* \mid \delta^*(q, w) \in A\}$.
Two states $q$ and $r$ are \emph{right language equivalent} if $RL(q) = RL(r)$.
If $RL(q) \neq RL(r)$ then one says $q$ and $r$ are \emph{distinguished}.
The language recognized by the automaton $(i, A, \delta)$ is the right language
$RL(i)$.
\paragraph{Acyclicity}
An DFA is \emph{acyclic} (an ADFA) if there is no non-empty sequence of
transitions from a state to it self. Formally, one has $\delta^*(q, w) =
q$ \textit{if and only if} $w = \epsilon$, for any state $q \in N$.
\paragraph{Reachability} A DFA is \emph{reachable} if any state is
\emph{reachable} from the \emph{initial state}. That means there exists a word
$w \in \Sigma^*$ such that $\delta^*(i, w) = q$ for any state $q$.
\paragraph{Coreachability} A DFA is \emph{coreachable} if for all state reachs
an \emph{accepting state}. That means there exists $w$ such that $\delta^*(q,
w) \in A$ for any state $q$.
\paragraph{Non initial DFA} We extend the definition of DFA to structures 
$(A, \delta)$ without \emph{initial state}.
\paragraph{Minimal DFA}
A DFA is \emph{minimal} if there is no DFA with fewer state which recognizes the
same language.\\
\setlength\parskip{0pt}

%%%%%%%%%%%%%%%
\subsection{Minimal ADFA}
%%%%%%%%%%%%%%%
An important point will be the notion of \emph{simple DFA}:
\begin{defi}
    A DFA is \emph{simple} if all its states are \emph{distinguished}.
\end{defi} 
\begin{prop}
    If a DFA is \emph{simple} then it doesn't have a non-trivial automorphism.
\end{prop}
This proposition expresses the problem of counting M(A)DFA having labeled or
unlabeled states are equivalent. So in the following, we consider automata
as labeled combinatorial structures/objects. Moreover, from the definition of
simple automata, it is easy to use the \textsc{Myhill-Nerode} theorem about 
\emph{minimal DFA}:
\begin{thm}[Myhill-Nerode]
    \label{thm:myhill-nerode}
    A DFA is \emph{minimal} if and only if it is \emph{reachable},
    \emph{coreachable}, and \emph{simple}.
\end{thm}
In the following, we enumerate \emph{non-initial ADFA} which are
\emph{coreachable} and \emph{simple}. Following \cite{liskovets2006exact}, we
fix an initial state and extract the underlying \emph{reachable} ADFA.
Ultimately, the extracted ADFA is \emph{reachable}, \emph{coreachable} and
\emph{simple}, and so it is minimal.
%%%%%%%%%%%%%%%
\subsection{Non-initial ADFA}
%%%%%%%%%%%%%%%
%
\begin{prop}
    Let $\aut$ be an ADFA.
    
    The \emph{absorbing state} $\emptyset$ is always \emph{reachable} from any
    state of $\aut$.
\end{prop}
This proposition is obvious because there is no loop in an ADFA.
For the upcoming bijection with parking functions, this point will be
important.
In particular, the bijection is based on the fact that there exists an
order on states wherein $\emptyset$ is minimum. From this order we
will put forward the states $q$ such that $\delta(q,a) = \emptyset$ for any
symbol $a \in \Sigma$ and which are \emph{accepting} or not.
\begin{lemme}
    The automaton $\aut$ is \emph{coreachable} if all state $q$ such that for
    any $a \in \Sigma$, one has $\delta(q,a) = \emptyset$ are \emph{accepting}.
\end{lemme}
The interpretation of states in parking functions will give a simple
caracterization of generalized parking functions associated to (non-initial)
\emph{coreachable ADFA}. Furthermore the flexibility of the generalized
parking function definition will easily give a family of parking functions,
exactly encoding those coreachable ADFA.

We denote $\niADFA^k_n$ the set of non-initial ADFA with $n$ labeled states (and
one absorbing state) over an alphabet of $k$ symbols and call it the graded
component of degree $n$.
We also denote $\niADFA^k := \bigsqcup_{n \seq 1} \niADFA^k_n$ the graded set 
of all \emph{non-initial ADFA}. Likewise, we
denote $\niscADFA^k$ the graded set of \emph{non-initial simple coreachable
ADFA}. This first set $\niADFA^k$ will be usefull to recall the extended notion
of non-initial ADFA given in \cite[Quasi-acyclic automata, \SOld
2.15]{liskovets2006exact} and the second to define constraints on extended
non-initial ``simple coreachable'' ADFA.
Considering a non-initial ADFA $\aut$ and a fixed state $i$, the sub-automata
$\aut^{(i)}$ extracted from all reachables states from $i$ defines an
(initial-connected) ADFA. We will extend ADFA to caracterize the complement
of $\aut^{(i)}$. Furthermore, an important remark is that the automaton
$\aut^{(i)}$ is \emph{minimal} if $\aut$ is \emph{coreachable} and
\emph{simple} (Theorem \ref{thm:myhill-nerode}):
\begin{align}
    \vcenter{\hbox{\scalebox{.8}{\newcommand{\AOnooo}{\node (0) [circle, draw,
    double] {$1$} ;}\newcommand{\AOnooa}{\node (1) [circle, draw, double] {$3$}
;}\newcommand{\AOnoob}{\node (2) [circle, draw, double] {$2$}
;}\newcommand{\AOnooc}{\node (3) [circle, draw, double] {$4$}
;}\newcommand{\AOnood}{\node (4) [circle, draw, fill=gray,text=white] {$5$}
;}\newcommand{\AOnooe}{\node (5) [circle, draw] {$7$}
;}\newcommand{\AOnoof}{\node (6) [circle, draw] {$8$}
;}\newcommand{\AOnoog}{\node (7) [circle, draw] {$6$}
;}\newcommand{\AOnooh}{\node (8) {$\emptyset$}
;}\begin{tikzpicture}[auto]
\matrix[column sep=1cm, row sep=.5cm, ampersand replacement=\&]{
 \AOnoob  \&\AOnooe \& \AOnooo \&  \AOnooh \\ 
  \AOnoog \&        \& \AOnooc \&  \\ 
 \AOnood  \&\AOnooa \&         \& \AOnoof \\ 
};
\path[->] 
     (2) edge[bend left=40] node[swap] {$\scriptstyle b$} (0)
     (6) edge node[swap] {$\scriptstyle a$} (8)
     (2) edge[bend left=40] node {$\scriptstyle a$} (8)
     (7) edge node {$\scriptstyle a$} (3)
     (4) edge node[swap] {$\scriptstyle b$} (7)
     (7) edge node[swap] {$\scriptstyle b$} (2)
     (6) edge node {$\scriptstyle b$} (3)
     (0) edge node {$\scriptstyle b$} (8)
     (1) edge node {$\scriptstyle a$} (6)
     (4) edge[bend right=40] node {$\scriptstyle a$} (6)
     (5) edge node {$\scriptstyle a$} (7)
     (1) edge node[swap] {$\scriptstyle b$} (7)
     (0) edge node[swap] {$\scriptstyle a$} (3)
     (5) edge node {$\scriptstyle b$} (0)
     (3) edge node[swap] {$\scriptstyle a,b$} (8);
\end{tikzpicture}}}} \quad\rightsquigarrow\quad 
\vcenter{\hbox{\scalebox{.8}{\newcommand{\AOnooo}{\node (0) [circle, draw, initial] {$5$}
;}\newcommand{\AOnooa}{\node (1) [circle, draw] {$6$}
;}\newcommand{\AOnoob}{\node (2) [circle, draw] {$8$}
;}\newcommand{\AOnooc}{\node (3) [circle, draw, double] {$4$}
;}\newcommand{\AOnood}{\node (4) [circle, draw, double] {$2$}
;}\newcommand{\AOnooe}{\node (5) [circle, draw, double] {$1$}
;}\newcommand{\AOnoof}{\node (6) {$\emptyset$}
;}\begin{tikzpicture}[auto, initial text=]
\matrix[column sep=1cm, row sep=.5cm,ampersand replacement=\&]{
 \AOnooo \& \AOnooa \& \AOnood \& \AOnooe \\ 
         \& \AOnoob \& \AOnooc \& \AOnoof \\ 
};
\path[->] %(4) edge node {$\scriptstyle a$} (6)
     (5) edge node[swap] {$\scriptstyle a$} (3)
     (1) edge node {$\scriptstyle b$} (4)
     (4) edge node {$\scriptstyle b$} (5)
     (1) edge node {$\scriptstyle a$} (3)
     (2) edge[bend right=40] node {$\scriptstyle a$} (6)
     (2) edge node {$\scriptstyle b$} (3)
     (3) edge node {$\scriptstyle a,b$} (6)
     (0) edge node {$\scriptstyle b$} (1)
     (5) edge node {$\scriptstyle b$} (6)
     (0) edge node {$\scriptstyle a$} (2);
\draw[->] (4) .. controls (3.5,1.5) and (4,1.5) .. node {$\scriptstyle a$}(6);
\end{tikzpicture}}}}
    \label{eq:extrac_A_i_of_A}
\end{align}
A good definition of extended (coreachable and simple) ADFA (with constraint)
gives a bijection between non-initial (coreachable and simple) ADFA and couple
of extended (coreachable and simple) ADFA (with constraint) and (M)ADFA.
%%%%%%%%%%%%%%%
\subsection{Extended non-initial ADFA}
%%%%%%%%%%%%%%%
In \cite{liskovets2006exact}, the author introduces ``quasi-acyclic automata''.
We call these objects extended non-initial ADFA (with $t$ extra
absorbing states), which means one considers tuples $(A, T, \delta)$ with $A$ as
the \emph{accepting state set}, $T$ the set of extra absorbing states, and
$\delta$ an extension of the transition function definition:
\begin{align*}
    \delta : N \times \Sigma \quad \longrightarrow \quad N \cup \{\emptyset\}
    \cup T
\end{align*} 
with $N$ the set of $n$ (labeled) states, $\Sigma$ the alphabet. 
One denotes by $\eADFA^{k,t} := \bigsqcup_{n \seq 0} \eADFA^{k,t}_n$ the graded
set of \emph{extended non-initial ADFA} (with $t$ extra absorbing states over
an alphabet of $k$ symbols).\\[5pt]
%\begin{ex}\ \\
\begin{minipage}[t!]{.55\textwidth}
\begin{ex}
 This example represents an extended non-initial ADFA with $3$ extra
 absorbing states $T = \{\alpha_1, \alpha_2, \alpha_3\}$ over the alphabet
 $\{a,b\}$.
 This structure $(\{2,4\}, \delta)$ is in $\eADFA^{2,3}_5$ with $\delta(1,a) =
 2$, $\delta(1,b) = \alpha_1$, $\delta(2, a) = \delta(2,b) = 4$, and so on.
 \end{ex}
\end{minipage}
\hfill
\begin{minipage}[t!]{.4\textwidth}
\scalebox{.8}{\newcommand{\AOnooo}{\node (0) [circle, draw] {$1$}
;}\newcommand{\AOnooa}{\node (1) [circle, draw,double] {$2$}
;}\newcommand{\AOnoob}{\node (2) [circle, draw] {$3$}
;}\newcommand{\AOnooc}{\node (3) [rectangle, draw, dashed] {$\alpha_1$}
;}\newcommand{\AOnood}{\node (4) [circle, draw,double] {$4$}
;}\newcommand{\AOnooe}{\node (5) [circle, draw] {$5$}
;}\newcommand{\AOnoof}{\node (6) {$\emptyset$}
;}\newcommand{\AOnoog}{\node (7) [rectangle, draw, dashed] {$\alpha_2$}
;}\newcommand{\AOnooh}{\node (8) [rectangle, draw, dashed] {$\alpha_3$}
;}\begin{tikzpicture}[auto]
\matrix[column sep=1cm, row sep=.5cm,ampersand replacement=\&]{
 \AOnooo \&         \& \AOnooc \& \AOnooe \& \AOnooh \\ 
 \AOnoob \& \AOnooa \& \AOnood \& \AOnoog \& \AOnoof \\ 
};
\path[->] (4) edge node {$\scriptstyle a$} (5)
     (0) edge node {$\scriptstyle b$} (3)
     (0) edge node {$\scriptstyle a$} (1)
     (2) edge node {$\scriptstyle a$} (1)
     (1) edge node {$\scriptstyle a,b$} (4)
     (2) edge[bend right] node {$\scriptstyle b$} (7)
     (4) edge node {$\scriptstyle b$} (3)
     (5) edge node {$\scriptstyle a$} (6)
     (5) edge node {$\scriptstyle b$} (8)
     (4) edge node {$\scriptstyle b$} (7);
\end{tikzpicture}}
\end{minipage}
%\end{ex}
\begin{remarq}$\!\!$
    An non-initial ADFA is an extended non-initial ADFA with $0$ extra absorbing
    states, %so one has 
    \mbox{$\niADFA^k \!= \eADFA^{k,0}$}.
\end{remarq}
%%%%%%%%%%%%%%%%%
\subsubsection{Enumeration of underlying transition functions}
%%%%%%%%%%%%%%%%%
In \cite{liskovets2006exact}, the author gives a formula $\dadfa$ to enumerate
the number of extended transition function $\delta$, underlying an extended
non-initial ADFA:
\begin{align*}
    \dadfa(k,t;\ n) = \sum_{j=0}^{n-1} \binom{n}{j} (-1)^{n-j-1} (j +
    t + 1)^{k(n-j)} \dadfa(k,t;\ j)\,, \tag*{\cite[Theorem
    3.1]{liskovets2006exact}}
\end{align*}
with $k$ the cardinal of the alphabet, $t$ the number of extra absorbing states
and $n$ the number of states.
In the same way this formula can be adapted to enumerate $\eADFA^{k,t}_n$:
\begin{coro}[of \expandafter{\cite[Theorem 3.1]{liskovets2006exact}}]
    \label{coro:enum_extended_non_initial_ADFA}
    The extended non-initial ADFA with $n$ states and $t$ extra absorbing
    states over an alphabet of $k$ symbols, $\eADFA^{k,t}_n$, are enumerated by
    the formula:
\begin{align*}
    \eadfa(k,t;\ n) = \sum_{j=0}^{n-1} \binom{n}{j} (-1)^{n-j-1} (2(j +
    t + 1)^k)^{n-j} \eadfa(k,t;\ j)\,,
\end{align*}
    for any $n \seq 1$ and $\eadfa(k,t;\ 0) = 1$.
\end{coro}
\begin{preuve}
    Immediate from the fact $\eadfa(k,t,n) = 2^n\dadfa(k,t,n)$.
\end{preuve}
In the next section we use this formula to show their there is an isomorphism
with some generalized parking functions using formula in
\cite{kung2003goncarov} (also recall in \cite{priezvirmaux}).
%
%
%%%%%%%%%%%%%%%%%
\subsubsection{Enumeration of ADFA}
%%%%%%%%%%%%%%%%%
In this subsection, we recall a well-know method of counting connected graphs,
used in \cite[Theorem 3.2]{liskovets2006exact}. This method points out that, for
a fixed state $i$, for any \emph{non-initial ADFA} $\aut$, one has a
reversible splitting of $\aut$ into an \emph{ADFA} $\aut^{(i)}$ (initially
connected by a fixed state $i$, see (\ref{eq:extrac_A_i_of_A})), and its
complement $\widebar{\aut}^{(i)}$ which is an extended non-initial ADFA.
It results an enumeration formula $\adfa(k; t)$ of ADFA over an alphabet of $k$
symbols with $t$ states (and one fixed label: the initial state). This formula
is given by the linear recurrence:
\begin{align*}
    \eadfa(k,0;\ n) = \sum_{t = 1}^{n} \binom{n-1}{t-1}\eadfa(k, t;\ n-t)
    \adfa(k;\ t)\,. \tag*{\cite[Theorem 3.2]{liskovets2006exact}}
\end{align*}
The complement $\widebar{\aut}^{(i)}$ is defined as follow:
let $\aut\in \niADFA^k_n$ be a non-initial ADFA with state set is $N$ and
let $\aut^{(i)} = (i, A_i, \delta_i)$ be an ADFA with its state set
is $N_i$ (the reachable states from $i$) and $\delta_i$ is the restriction of
$\delta$ to states of $N_i$.
We set $\widebar{\aut}^{(i)}$ the complement of $\aut^{(i)}$ as an
\emph{extended non-initial ADFA} of $n-t$ states $\widebar{N}_i = N \backslash
N_i$ with $t$ extra absorbing states $N_i$, the accepting states are
$\widebar{A}_i = A \cap \widebar{N}_i$ and the extended transition function
$\widebar{\delta}_i$ is defined by:
\begin{align*}
    \widebar{\delta}_i : \begin{array}{rcl}
        \widebar{N}_i \times \Sigma & \longrightarrow & N \cup
        \{\emptyset\}\,,\\
        q \times a & \longmapsto & \delta(q,a)\,.
    \end{array}
\end{align*}
We denote $\split_i : \niADFA^k \to \ADFA^k \times \eADFA^{k}$ that
splitting bijection (with $\eADFA^k = \bigsqcup_{t \seq 1} \eADFA^{k,t}$). The
inverse bijection consists simply (in term of graph) to merge extra absorbing
states $q$ of $\widebar{\aut}^{(i)}$ with the state $q$ of $\aut^{(i)}$.
\begin{remarq}
    $N \cup \{\emptyset\}= \widebar{N}_i \cup \{\emptyset\} \cup N_i$.
\end{remarq}
\begin{minipage}[t!]{.55\textwidth}\color{fontcolor}
\begin{ex}
 In Equation (\ref{eq:extrac_A_i_of_A}), we have a non-initial ADFA $\aut$ on
 the left and $\aut^{(5)}$ on the right. We represent (here, on the right)
 $\widebar{\aut}^{(5)}$, the complement of $\aut^{(5)}$, as an extended ADFA
 (with extra absorbing states framed with dashed rectangles). 
\end{ex}
\end{minipage}
\hfill
\begin{minipage}[t!]{.4\textwidth}
    \begin{center}
% 5 6 2 8 1 4
% 4 5 2 6 1 3
    \scalebox{.8}{\newcommand{\AOnooo}{\node (0) [circle, draw, double] {$3$}
;}\newcommand{\AOnooa}{\node (1) [circle, draw] {$7$}
;}\newcommand{\AOnoob}{\node (2) [rectangle, draw, dashed] {$5$}
;}\newcommand{\AOnooc}{\node (3) [rectangle, draw, dashed] {$2$}
;}\newcommand{\AOnood}{\node (4) [rectangle, draw, dashed] {$1$}
;}\newcommand{\AOnooe}{\node (5) [rectangle, draw, dashed] {$4$}
;}\newcommand{\AOnoof}{\node (6) [rectangle, draw, dashed] {$6$}
;}\newcommand{\AOnoog}{\node (7) [rectangle, draw, dashed] {$8$}
;}\newcommand{\AOnooi}{\node (9) {$\emptyset$}
;}\begin{tikzpicture}[auto]
\matrix[column sep=1cm, row sep=.5cm, ampersand replacement=\&]{
 \AOnoob \& \AOnooc \& \AOnooe \& \AOnooi \\ 
 \AOnood \& \AOnooa \& \AOnoof \& \AOnooo \& \AOnoog \\ 
};
\path[->] (0) edge node {$\scriptstyle a$} (7)
              edge node[swap] {$\scriptstyle b$} (6)
          (1) edge node {$\scriptstyle a$} (6)
              edge node[swap] {$\scriptstyle b$} (4);
\end{tikzpicture}}
    \end{center}
\end{minipage}
%\end{ex}
%%%%%%%%%%%%%%%
\subsection{Extended coreachable simple non-initial ADFA with constraints}
    \label{ssec:extended_corea_simpl_non_initial_ADFA_constraints}
%%%%%%%%%%%%%%%
In this subsection, we focus on \emph{non-initial ADFA}
which are \emph{coreachable} and \emph{simple}. We start by giving a definition
of \emph{extended simple non-initial ADFA}:
\begin{defi}
    An \emph{extended non-initial ADFA} $(A, \delta, T)$ is \emph{simple} if 
    one has $RL(q) \neq RL(r)$ or there exists $w \in \Sigma^*$
    such that $\delta^*(q, w) \neq \delta^*(r, w)$ with $\delta^*(q, w) \in T$,
    for any distincts states $q, r$.
\end{defi}
This definition is another way to say that if $\split_i$ is applied on simple
ADFA $\aut$ then it gives a couple of simple structures: $(\aut^{(i)},
\widebar{\aut}^{(i)})$.
Unlike the splitting of \cite{liskovets2006exact} described in the previous
subsection, the restriction of $\split_i$ to $\niscADFA^k$ is not a bijection.
From any couple of simple structures an extended coreachable simple
non-initial ADFA and a MADFA do not necessary produce a simple non-initial
ADFA.
The reason is (without well-chosen constraints) for some MADFA $\aut$ and some
extended simple (and coreachable) ADFA $\autb$ could share the ``same''
transitions on their respective states. Due to this therefore we add
constraints to obtain a reversible splitting.
%%%%%%%%%%%%%%%
%%%%%%%%%%%%%%%
\subsubsection{Preservation of simplicity and constraints}
%%%%%%%%%%%%%%%
To be sure that the assembly of a couple of an \emph{extended simple non-initial
ADFA} and a \emph{simple ADFA} remains \emph{simple} (non-initial ADFA), one
considers the couple of extended ADFA, and ADFA satisfying a set of constraints
$C$.
The idea is to forbid them from sharing states with the same transitions, and
the same accepting status.
A set of constraints $C$ is a set of couples $(\nu, b)$ with $\nu: \Sigma \to T
\cup \{\emptyset\}$ and $b$ an accepting status (\textit{true} or \textit{false}). A
couple $(\aut, \autb)$ of ADFA and extended ADFA satisfies $C$ if
\begin{itemize}
  \item The state set of $\aut$ is $T$ and for each state $q$ of $\aut$ there
    exists an unique couple $(\nu, b)$ such that $\delta_q = \nu$ and $q$ is
    accepting if $b$ is \textit{true},
  \item The set of extra absorbing states of $\autb$ is $T$ and for any
  state $q$ of $\autb$ and any couple $(\nu, b)$ of $C$, one has $\delta_q \neq
  \nu$ or the accepting status of $q$ is the negation of $b$.
\end{itemize}
\begin{prop}
    For any couple $(\aut, \autb)$ of simple structures satisfying $C$, one has
    $\split_i^{-1}(\aut, \autb)$ is simple.
\end{prop}
\begin{remarq}
    There does not always exist a couple of ADFA and extended ADFA that
    satisfies any set of constraints.
\end{remarq}
%%%%%%%%%%%%%%%
%%%%%%%%%%%%%%%
\subsubsection{Preservation of coreachability}
%%%%%%%%%%%%%%%
Both structures $\aut$ and $\autb$ must not contain state $q$ such that
$\delta(q, a) = \emptyset$ (for any symbol $a$) and $q$ is not accepting.
Furthermore $\autb$ must not contain state $q$ such that $\delta(q, a) =
\emptyset$ at all. Otherwise, $\split^{-1}_i(\aut, \autb)$ is not
coreachable.
So we define the notion of \emph{coreachability} of extended ADFA by:
\begin{defi}
    An extended ADFA is \emph{coreachable} if there is no state $q$ such that
    $\delta(q,a) = \emptyset$, for any $a \in \Sigma$. 
\end{defi}
We denote by $\esADFA^{k,C}$ and $\mADFA^{k, C}$ respectively the graded set of
\emph{extended simple coreachable ADFA} and of \emph{MADFA} satisfying $C$.
Due to this is a well-chosen definition of extended coreachable and simple
ADFA, we obtain the bijection:
% one obtains an isomorphism:
%
\begin{lemme}
    \label{lem:niscADFA_simeq_MADFA_esADFA}%\ \\[-30pt]
    \begin{align*}
        \niscADFA^k \simeq \bigsqcup_{C} \mADFA^{k, C} \times \esADFA^{k, C}\,.
    \end{align*}
\end{lemme}
This gives us a description of coreachable and simple non-initial ADFA as the
direct sum (over all sets of constraints $C$) of couple MADFA, and extended
simple coreachable ADFA satisfying $C$. In the next section, the set of
all available constraints will become aparent.

%% file: sec_pf.tex
In this section, we recall the constructive definition of generalized parking
functions given by \textsc{Virmaux} and the author in \cite{priezvirmaux}. In a
first part we define a bijection between of \emph{non-initial ADFA} and a
remarkable family of \emph{generalized parking functions}.
This bijection will reveal two interesting points:
\begin{itemize}
  \item The localization of (some) \textit{non-distinguishable} states.
  According to this, one extracts a sub-family of parking functions that are isomorphic to the
      \emph{simple non-initial ADFA}.
  \item An easy translation of some constraints. In particular, in this first
      part, we give an analoguous family of parking functions which are
      isomorphic to the \emph{coreachable non-initial ADFA}.
\end{itemize}
By combining both  of these points we obtain an isomorphism with
\emph{(extended) simple coreachable ADFA}. Finally in the next section we go
back to Lemma \ref{lem:niscADFA_simeq_MADFA_esADFA} to formulate a recurrence
relation which enumerate MADFA.
 %
%%%%%%%%%%%%%%%
\subsection{Definition}
%%%%%%%%%%%%%%%
Parking functions were first introduced in \cite{konheim1966occupancy} to
model hashing problems in computer science and appear in many different
contexts in combinatorics, such as \textit{labeled trees}, \textit{prüfer
sequence}, \textit{hyperplane arrangements}, \textit{etc}. A parking function on a finite
set $N$ is a function $f : N \to \NN_+$ such that $\card{f^{-1}([k])} \seq k$,
for any $k \in [n]$ (with $n = \card{N}$ and $[k] = \{1, \cdots, k\}$).
A generalization of parking functions was formulated in
\cite{stanley2002polytope} and well studied in \cite{kung2003goncarov,
priezvirmaux}. Let $\seqPF : \NN_+ \to \NN$ be an non-decreasing function; a
$\seqPF$-parking function is a function $f$ such that
$\card{f^{-1}([\seqPF(k)])} \seq k$, for any $k \in [n]$.
\begin{remarq}
    Usual parking functions are $\seqPF$-parking function with $\seqPFn{k} =
    k$.
\end{remarq}
% We fix a convention that $\seqPF$ could be extend as a map from $\NN$ into $\NN$
% with $\seqPFn{0} = 0$.
%%%%%%%%%%%%%%%%%%%%%%%%%%%%%%%%%%%%
\subsection{Constructive definition}
%%%%%%%%%%%%%%%%%%%%%%%%%%%%%%%%%%%%
In this paper, we encode $\seqPF$-parking functions by the sequence $(Q_j)$
defined by $Q_j = f^{-1}(\{j\})$. We therefore define $\seqPF$-parking functions
on $N$ as a sequence of $\seqPFn{n}$ disjoint subsets $(Q_j)$ of $N$ satisfying:
\begin{align}
    \sum_{i = 1}^{\seqPFn{k}} \card{Q_i} \seq k\,,\quad\qquad \text{for any}\; k
    \in [n]\,.
    \label{eq:pf_condition}
\end{align}
\begin{remarq}
    The parking function condition imposes that $\card{f^{-1}(\{k\})} = 0$ for
    any $k > \seqPFn{n}$. So the definition in terms of a set sequence allows
    to complete the sequence with an arbitrary sequence of empty sets. However
    one considers $\seqPF$-parking functions as finite sequences of sets.
\end{remarq}
The main advantage of this definition (in terms of sequences of sets) is that it
involves a natural recursive definition (see
\cite[\pvsectiondefrec]{priezvirmaux}).
A convenient language for this is the \emph{species theory}
\cite{bergeron1998combinatorial} (or equivalently \emph{decomposable
combinatorial class} \cite{flajolet2009analytic}). 
Let $\speciesSet$ be the species of sets (such that $\speciesSet[U] := \{U\}$
for any finite set $U$), let $\speciesOne$ be the characteristic empty species
(such that $\speciesOne[U]  = \{\emptyset\}$ if $U = \emptyset$ and
$\emptyset$ in otherwise).
We denote by $+$ (and $\sum$) the sum of species (disjoint union of labeled
combinatorial structures: $(\speciesA + \speciesB)[U] = \speciesA[U] \cup
\speciesB[U]$), by $\cdot$ (, $\prod$ and the exponentiation) the product of
species (cartesian product of labeled combinatorial structures: $(\speciesA
\cdot \speciesB)[U] = \sum_{S \sqcup T = U} \speciesA[S] \times \speciesB[T]$).

We directly give the recursive solution (of \cite[Eq.
\pveqrecursiveconstruction]{priezvirmaux}) which defines $\seqPF$-parking
functions, $\PF(\seqPF)$, grade by grade (see
\cite[\pvcorsolutionrecurrence]{priezvirmaux}) as a sum over all compositions
$\pi$ of the integer $n$, noted $\pi \models n$:
\begin{align}
    \PF(\seqPF) = \speciesOne + \sum_{n \seq 1} \PF_n(\seqPF)&&
    \text{with}\qquad \PF_n(\seqPF) = \sum_{\pi \models n}
    \prod_{i=1}^{\ell(\pi)} \left(\speciesSet^{\Upsilon(\seqPF; \pi,
    i)}\right)_{\pi_i}\,,
    %\tag*{\cite[\pvcorsolutionrecurrence]{priezvirmaux}}
    \label{eq:pf_constructive_def}
\end{align} 
where $\Upsilon$ is
\begin{align}
    \Upsilon(\seqPF; \pi, i) := \begin{dcases*}
        \seqPFn{1} & if $i = 1$,\\
        \seqPFn{1 + \pi(i-1)} - \seqPFn{1 + \pi(i-2)} & otherwise,
    \end{dcases*}
    \label{eq:upsilon_def} 
\end{align} 
with $\pi(i) = \pi_1 + \cdots + \pi_i$ the partial sum of the first $i$ parts of
$\pi$.
\begin{ex}
    \label{ex:expansion_PF}
    Here is the constructive definition of $\PF_3(m^2)$ expanding
 (\ref{eq:pf_constructive_def}):
 \begin{align*}
    % 3, 21, 12, 111
    \PF_3(m^2) = \left(\speciesSet^{1}\right)_{3} +
    \left(\speciesSet^{1}\right)_{2}\cdot \left(\speciesSet^{8}\right)_{1} +
    \left(\speciesSet^{1}\right)_{1}\cdot \left(\speciesSet^{3}\right)_{2} +
\left(\speciesSet^{1}\right)_{1}\cdot \left(\speciesSet^{3}\right)_{1} \cdot
\left(\speciesSet^{5}\right)_{1}
 \end{align*}
\end{ex}
By abuse of notation we denote $\seqPF$ directly by its image over $m$. 
For example, we denote $\PF(m^2)$ the generalized parking functions
$\PF(\seqPF)$ with $\seqPFn{m} := m^2$.
Furthermore by abuse again, we identify $\PF_n(\seqPF)$ to the
$\PF(\seqPF)$-structures on the set $[n]$.
\begin{ex}
    \label{ex:pf_structures_012_3}
 We denote $(\{a,b,c, \cdots\}, \{d,e,f, \cdots\}, \ldots)$ by
    $\boldsymbol{(}abc\cdots \mid def\cdots \mid \ldots\boldsymbol{)}$. 
 The first $\PF_n(m^2)$-structures for $n=0, 1$ and $2$ are:\\ %[-10pt]
\begin{minipage}[t!]{.3\textwidth}
\begin{align*}
    \PF_0(n^2) = \{\;& () \;\}\,,\\
    \PF_1(n^2) = \{\;& (1) \;\}\,,
\end{align*}
\end{minipage}
\hfill
\begin{minipage}[t!]{.65\textwidth}
\begin{align*}
    \PF_2(n^2) = \{\;& 
      \boldsymbol{(}12\mid \cdot\mid\cdot\mid\cdot\boldsymbol{)},
      \boldsymbol{(}1\mid 2\mid\cdot\mid\cdot\boldsymbol{)},
      \boldsymbol{(}1\mid\cdot\mid 2\mid\cdot\boldsymbol{)},
      \boldsymbol{(}1\mid\cdot\mid\cdot\mid2\boldsymbol{)},\\&
      \boldsymbol{(}2\mid1\mid\cdot\mid\cdot\boldsymbol{)},
      \boldsymbol{(}2\mid\cdot\mid1\mid\cdot\boldsymbol{)},
      \boldsymbol{(}2\mid\cdot\mid\cdot\mid1\boldsymbol{)}
    \;\}\,.
\end{align*}
\end{minipage}
And from Example \ref{ex:expansion_PF}, some of the $27$ structures of
$\PF_3(m^2)$ resulting from $(\speciesSet^1)_2\cdot (\speciesSet^8)_1$ are:
\begin{align*}
\boldsymbol{(}12\mid 3\mid\cdot\mid\cdot\mid\cdot\mid\cdot\mid\cdot\mid\cdot\mid\cdot\boldsymbol{)},\  
\boldsymbol{(}12\mid \cdot\mid 3\mid\cdot\mid\cdot\mid\cdot\mid\cdot\mid\cdot\mid\cdot\boldsymbol{)},\ 
\cdots,\ 
\boldsymbol{(}12\mid \cdot\mid\cdot\mid\cdot\mid\cdot\mid\cdot\mid\cdot\mid\cdot\mid 3\boldsymbol{)},\\
\boldsymbol{(}13\mid 2\mid\cdot\mid\cdot\mid\cdot\mid\cdot\mid\cdot\mid\cdot\mid\cdot\boldsymbol{)},\ 
\boldsymbol{(}13\mid \cdot\mid2\mid\cdot\mid\cdot\mid\cdot\mid\cdot\mid\cdot\mid\cdot\boldsymbol{)},\  
\cdots,\ 
\boldsymbol{(}23\mid \cdot\mid\cdot\mid\cdot\mid\cdot\mid\cdot\mid\cdot\mid \cdot\mid1\boldsymbol{)}.  
\end{align*}
\end{ex}
%
%%%%%%%%%%%%%%%%%%%%%%%%%%%%%%%%%%%%
\subsection{Enumeration and interpretation} \label{ssec:enum_interpretation}
%%%%%%%%%%%%%%%%%%%%%%%%%%%%%%%%%%%%
In \cite{kung2003goncarov} the authors gave a recurrence relation to
enumerate $\PF_n(\seqPF)$-structures:
\begin{align*}
    \pf(\seqPF; n) = \sum_{j=1}^{n} (-1)^{j-1} \binom{n}{j} \seqPFn{n-j+1}^j
    \pf(\seqPF; n-j)\,.
    \tag*{\cite[Theorem 4.2]{kung2003goncarov}}
\end{align*}
From this formula and the Corollary \ref{coro:enum_extended_non_initial_ADFA},
we immediately obtain:
\begin{thm}
    \label{thm:iso_PF_ADFA}
    There is a bijection between $2(m+t)^k$-parking functions and
    \emph{extended non-initial ADFA} $\eADFA^{k,t}$.
\end{thm}
Thanks to \cite[Theorem 4.2]{kung2003goncarov} and Corollary
\ref{coro:enum_extended_non_initial_ADFA}.
% \begin{preuve}
%     Thanks to , one has:
%     \begin{align*}
%     \pf(2(m+t)^k; n) = \sum_{j=1}^{n} (-1)^{j-1} \binom{n}{j}
%     (2(n-j+1+t)^k)^j \pf(2(m+t)^k; n-j) = \eadfa(k,t;n)\,.
%     \end{align*}
% \end{preuve}
In the following we explicit the bijection. To do that we need to use a more
precise formula to extract simple ADFA/parking functions. 
In \cite{priezvirmaux}, we use the constructive definition
(\ref{eq:pf_constructive_def}) to obtain (automatically) this more expressive
formula, \textit{i.e.} the non-commutative Frobenius characteristic of the
natural action of the $0$-Hecke algebra on generalized parking functions
\cite[\pvthmFrobeniuscharacteristic]{priezvirmaux} (expressed in complete
$(\S^\pi)$-basis of non-commutative symmetric functions):
\begin{align}
    \ncch(\PF_n(\seqPF)) = \sum_{\pi \models n} \left(
        \sum_{\tau \models \ell(\pi)} \prod_{i=1}^{\ell(\tau)}
        \binom{\Psi_\tau(\seqPF; \pi, i)}{\tau_i}
    \right) \S^\pi
    \label{eq:ncch_PF}
\end{align}
with $\Psi_\tau$ a generalization of $\Upsilon$ (\ref{eq:upsilon_def}):
\begin{align*}
\Psi_\tau(\seqPF; \pi, i) = \begin{dcases*}
        \seqPFn{1} & if $i = 1$,\\
        \seqPFn{1 + \pi(\tau(i))} - \seqPFn{1 + \pi(\tau(i-1))} & in otherwise,
    \end{dcases*}
\end{align*}
\begin{remarq}
    \label{rq:interpretation_complete_ncsf}
    The complete non-commutative symmetric functions $(\S^\pi)$ are a
    convenient algebraic way to encode the action of relabeling of set
    sequence. The coefficient of $\S^\pi$ (with $\pi$ a composition of
    $n$) is the number of $\seqPF$-parking functions of $n$ (upto isomorphism) 
    such that the first non-empty set contains $\pi_1$ elements, the second one contains 
    $\pi_2$ elements, and so forth.
\end{remarq}
\begin{ex}
%     The complete function associated to $(\speciesSet^1)_3$ is $\S^3$, the
%     terms of $(\speciesSet^1)_2 \cdot (\speciesSet^8)_1$ is $8\S^{21}$ (with
%     $8$ is all distincts ways to dispatch one element in $8$ sets) and one has
%     $3\S^{12} + 3\S^{111}$ for $(\speciesSet^1)_1\cdot (\speciesSet^3)_2$ with
%     the first term $\S^{12}$ encodes sets with one elements in the first part
%     and two dispatch in the same set between the last three sets, the second
%     terms encodes that one has dispatch two elements in distincts sets between
%     the last three sets (see Examples \ref{ex:expansion_PF} and
%     \ref{ex:pf_structures_012_3}).
\begin{align*}
    % 3, 21, 12, 111
    \PF_3(m^2) &= \underbrace{\left(\speciesSet^{1}\right)_{3}}_{} +
    \underbrace{\left(\speciesSet^{1}\right)_{2}\cdot \left(\speciesSet^{8}\right)_{1}}_{} +
    \underbrace{\left(\speciesSet^{1}\right)_{1}\cdot \left(\speciesSet^{3}\right)_{2}}_{} +
\underbrace{\left(\speciesSet^{1}\right)_{1}\cdot \left(\speciesSet^{3}\right)_{1} \cdot
\left(\speciesSet^{5}\right)_{1}}_{}\\[-8pt]
\ncch(\PF_3(m^2)) &= \quad\S^3\;\, + \qquad 8\S^{21}\quad\; +\;
3\S^{12} + 3\S^{111} + \qquad 15\S^{111}
 \end{align*}
\end{ex}
From the non-commutative characteristic (\ref{eq:ncch_PF}) the specialization
of $\S^\pi$ to the multinomial $\binom{n}{\pi_1, \cdots, \pi_k}$ gives
automatically another formula:
% \begin{coro}[of \expandafter{\cite[\pvthmFrobeniuscharacteristic]{priezvirmaux}}]
    \begin{align*}
    \pf(\seqPF; n) = \sum_{\pi \models n} \left(
        \sum_{\tau \models \ell(\pi)} \prod_{i=1}^{\ell(\tau)}
        \binom{\Psi_\tau(\seqPF; \pi, i)}{\tau_i}
    \right) \binom{n}{\pi_1, \cdots, \pi_{\ell(\pi)}}\,.
    \end{align*}
% \end{coro}
%\begin{remarq}
%     Unfortunately $\Psi_\tau$ is defined according to the part order of the two
%     compositions $\pi$ and $\tau$, so this new expression of $\pf$ is not
%     (computationally) efficient.
Unfortunately this new expression of $\pf$ is a double sum over compositions of
integers, so this is not (computationally) efficient.
%\end{remarq}
Finally we reuse that non-commutative characteristic
$\ncch(\PF_n(\seqPF))$ to extract an enumeration of a sub-family of
$\PF(\seqPF)$ according to \emph{simple ADFA} (\SOld \ref{ssec:simple_pf}).
Meanwhile, one explicits the bijection between $\PF(2(m+t)^k)$ and
$\eADFA^{k,t}$.
%%%%%%%%%%%%%%%%%%%%%%%%%%%%%%%%%%%%
\subsection{Explicit bijection $\PF(2(m+t)^k) \simeq \eADFA^{k,t}$}
%%%%%%%%%%%%%%%%%%%%%%%%%%%%%%%%%%%%
The following bijection is based on
the parking functions condition (\ref{eq:pf_condition}) involving 
a natural division of $(Q_j)$ into $n$ factors splitting between each
$\seqPFn{i}$.
This natural division will caracterize the sets of all transition functions
$\nu$ from an alphabet $\Sigma$ of $k$ symbols into a set of $p$ states fixing
one state $q$. This caracterisation will make sure that the
built automaton is acyclic.
%%%%%%%%%%%%%%%%%%%%%%%%%%%%
\subsubsection{Division factors} \label{ssec:division_factors}
%%%%%%%%%%%%%%%%%%%%%%%%%%%%
Let $N$ be a finite set of cardinal $n$.
%\begin{defi}
The \emph{division factors} of a $\seqPF$-parking function $(Q_j)$
on $N$ is defined as the sequence of factor $(D_p)_{p \in [n]}$ with $D_p =
(Q_j)$ with $j \in [\seqPFn{j-1}+1, \seqPFn{j}]$.
%\end{defi}
%
\begin{ex}
    \label{ex:division_factors}
    The $\PF(m^2)$-structures $\boldsymbol{(}\; 3 \mid \cdot \mid 1
    \mid \cdot \mid 2 \mid \cdot \mid \cdot \mid \cdot\mid \cdot
    \;\boldsymbol{)}$ has as \textit{division factors} $(D_1,D_2,D_3)$ with
    $D_1 = (\,3\,)$, $D_2 = (\,\cdot \mid 1 \mid \cdot\,)$ and $D_3 =(\,2 \mid \cdot
    \mid \cdot\mid \cdot \mid \cdot\,)$.
\end{ex}
%
%%%%%%%%%%%%%%%%%%%%%%%%%%%%
\subsubsection{Linear order on parking functions} \label{ssec:lin_order}
%%%%%%%%%%%%%%%%%%%%%%%%%%%%
Division factors will be used to fix an order between transitions according
to a ``fixing state''. From now on we set a total order on $N$ associated
to any $\PF(\seqPF)$-structures.
We suppose that we are given a total order $<_N$ on $N$
(in examples we use the natural order on $[n]$) and let $(Q_j)$ be a
$\seqPF$-parking function on $N$. We define a second total order $<_q$ by:
\begin{align*}
    q <_q q' \quad\Longleftrightarrow\quad \begin{dcases*}
        q \in Q_{k}\; \text{and}\; q' \in Q_{k'} & with $k < k'$ or\\
        q, q' \in Q_k & and $q <_N q'$\,.
    \end{dcases*} 
\end{align*}
\begin{remarq}
    The order $<_q$ is the linear order defined by the inverse of the
    standardization of the parking functions seen as words. (In Example
    \ref{ex:division_factors}, the parking function $(Q_j)$ could be represented
    by the word $w = 351$ (with $w_i = j$ \textit{iff} $i \in Q_j$) and the
    inverse of it standardization is $312$)
\end{remarq}
\begin{ex}
    \label{ex:order_pf_structures}
    The $m^2$-parking functions
    $\boldsymbol{(}\,
        4\mid 12\mid \cdot\mid \cdot\mid \cdot\mid \cdot\mid \cdot\mid \cdot\mid \cdot\mid 3\mid \cdot\mid \cdot\mid \cdot\mid \cdot\mid \cdot\mid \cdot
        \,\boldsymbol{)}$ on $[4]$
    and $\boldsymbol{(}\,
    3\mid\cdot\mid\cdot\mid25\mid\cdot\mid\cdot\mid\cdot\mid\cdot\mid\cdot\mid\cdot\mid\cdot\mid\cdot\mid\cdot\mid\cdot\mid1\mid\cdot\mid\cdot\mid\cdot\mid\cdot\mid\cdot\mid4\mid\cdot\mid\cdot\mid\cdot\mid\cdot
    \,\boldsymbol{)}$ on $[5]$
    define respectively the orders $4 <_q 1 <_q 2 <_q 3$
    and  $3 <_q 2 <_q 5 <_q 1 <_q 4$.
\end{ex}
%
%%%%%%%%%%%%%%%%%%%%%%%%%%%%
\subsubsection{Non-initial ADFA from $2m^k$-parking functions}
\label{ssec:adfa_from_pf}
%%%%%%%%%%%%%%%%%%%%%%%%%%%%
In this subsection, we associate to any $\PF(2m^k)$-structure a non-initial
ADFA $(A, \delta)$ of $\niADFA^k$. By definition this bijection will transport
interesting properties on the right language of states. This construction will
be easily generalized to (extended) non-initial ADFA with constraints, and
finally extended to coreachable simple non-initial ADFA.
Let $(Q_j)$ be a $\PF(2m^k)$-structures on $N$ and $(D_p)_{p \in [n]}$ its
\emph{division factors}. 
We complete the total order $<_q$ with an other element/state, the absorbing
state: $\emptyset <_q j$ for any $j \in N$. So we have $\emptyset
=: q_0 <_q q_1 <_q q_2 <_q \cdots <_q q_n$.
\begin{prop}
    There is exactly $p^k - (p-1)^k$ maps $\nu : \Sigma \to \{q_1, \cdots,
    q_{p-1}\} \cup \{\emptyset\}$ such that we have $\nu(a) = q_{p-1}$ for at
    least one $a \in \Sigma$.
\end{prop}
This proposition is obvious. By definition of $(D_p)$ each factor $D_p$
is a sequence of $2(p^k - (p-1)^k)$ sets. 
Each $D_p$ is associated to the set of all maps $\nu : \Sigma \to
\{q_1, \cdots, q_{p-1}\} \cup \{ \emptyset \}$ such that we have $\nu(a) =
q_{p-1}$ for at least one symbol $a \in \Sigma$.
(In others terms, we associate the unique map $\nu : \Sigma \to \{\emptyset\}$
to $D_1$, the maps $\nu : \Sigma \to \{q_1\}\cup \{\emptyset\}$ such that 
$\nu(a) = q_1$ for at least one $a \in \Sigma$, are associated to $D_2$ and so
on.)

For a fixed total order $<_{\Sigma}$ on $\Sigma$: $a_1 <_{\Sigma} a_2<_{\Sigma}
\cdots <_{\Sigma} a_k$, the lexicographical order on the
sequence of the image of $\nu$ seen as words: $\nu(a_1) \nu(a_2) \cdots
\nu(a_k)$ defines a total order on maps $\nu$: $\nu_1 < \nu_2 < \cdots < \nu_{p^k -
(p-1)^k}$ in each \textit{division factors}. We denote $\nu^{(p)}_j$ the maps
associated to the factor $D_p$.
We whizz up all that to finally define the bijection. Let $\bij$ be the map
which associates to any $(Q_j)$ the automaton structure $(A, \delta)$ defined
by:
\begin{itemize}
  \item the set of states is $N$ and the absorbing state is $\emptyset$,
  \item the set of accepting states $A$ is the union of $Q_j$ with $j$ even,
  \item the transition function $\delta$ is setting by $\delta_q :=
  \nu^{(p)}_j$ \textit{iff} $q \in Q_{2(p-1)^k+2j-1}$ or $q \in Q_{2(p-1)^k +
  2j}$.
\end{itemize}
%\begin{ex}\ \\
\begin{minipage}[t!]{.55\textwidth}
\begin{ex}
Let $(Q_i)$ be the $2m^2$-parking functions described by the division factors:
$D_1 = \boldsymbol{(}\;\cdot \;|\; 4\;\boldsymbol{)}$, 
$D_2= \boldsymbol{(}\; \cdot \;|\; \cdot \;|\; 
                     \cdot \;|\; \cdot \;|\;
                     3     \;|\; 1 \;\boldsymbol{)}$,
$D_3$ a sequence of $10$ empty sets and 
$D_4$ a sequence of $16$ sets with only the $7^\text{th}$ is non-empty set is:
$Q_{25} = \{2\}$ ($25 = (4-1)^2 + 7$).
\end{ex}
\end{minipage}
\hfill
\begin{minipage}[t!]{.4\textwidth}
\ \\The non-initial ADFA $\bij(Q_i)$.\\
$\vcenter{\hbox{\scalebox{.8}{\newcommand{\AOnooo}{\node (0) [circle, draw]
{$2$};}\newcommand{\AOnooa}{\node (1) [circle, draw] {$3$}
;}\newcommand{\AOnoob}{\node (2) [circle, draw, double] {$1$}
;}\newcommand{\AOnooc}{\node (3) [circle, draw, double] {$4$}
;}\newcommand{\AOnood}{\node (4)  {$\emptyset$}
;}\begin{tikzpicture}[auto]
\matrix[column sep=1cm, row sep=.1cm, ampersand replacement=\&]{
                \& \AOnooa \\ 
 \AOnooo \&                \& \AOnooc \& \AOnood \\ 
                \& \AOnoob \\ 
};
\path[->] 
     (1) edge node {$\scriptstyle a,b$} (3)
	 (3) edge node {$\scriptstyle a,b$} (4)
	 (0) edge node {$\scriptstyle b$} (2)
	 (0) edge node {$\scriptstyle a$} (1)
	 (2) edge node {$\scriptstyle a,b$} (3);
\end{tikzpicture}}}}$
\end{minipage}\\
% \end{ex}
Let $(Q_i)$ be a $2m^k$-parking function.
\begin{prop}
    The automaton structure $\bij(Q_i)$ is a non-initial ADFA of
    $\niADFA^{k}_n$.
\end{prop}
\begin{preuve}
    From the construction, each state has transitions defined, so it
    is a well defined non-initial DFA. Then the parking function condition
    (\ref{eq:pf_condition}) and the definition of $\delta$ assert that the
    structure is acyclic.
\end{preuve}
\begin{lemme}
    The map $\bij$ is a bijection.  
\end{lemme}
\begin{preuve}
    The different total orders ($<_N$, $<_q$, $<_\Sigma$ and the
    lexicographical order on maps $\nu$) assert that each non-initial ADFA is
    produced once and Theorem \ref{thm:iso_PF_ADFA} asserts that each is
    produced.
\end{preuve}
The map $\bij$ associates the same transitions to each state in the same set
$Q_j$. 
Furthermore, the order on $Q_j$ defined by the parking function involves
a ``height property'' as in \cite{revuz1992minimisation} to define an efficient
minimization. 
\begin{lemme}
    \label{lem:simple_niADFA_pf}
    The non-initial ADFA $\bij(Q_j)$ is simple and coreachable \textit{if and
    only if} $\card{Q_j} \ieq 1$, for any $j \in [2n^k]$ and $Q_1 = \emptyset$.
\end{lemme}
\begin{preuve}
    The idea is that if there exists distincts states $q$ and $s$ such that
    $RL(q) = RL(s)$ then there is a $Q_j$ such that $\card{Q_j} > 1$ in $(Q_j)$
    as in \cite{revuz1992minimisation}. So if $\card{Q_i} \ieq 1$ for any $i$ then
    $\bij(Q_i)$ is simple and the reciprocal is trivial. As seen in
    the previous section, a non-initial ADFA is coreachable if the states whose
    transitions going only to the absorbing state, are accepting.
\end{preuve}
%%%%%%%%%%%%%%%%%%%%%%%%%%%%
\subsubsection{Simple parking functions}
\label{ssec:simple_pf}
%%%%%%%%%%%%%%%%%%%%%%%%%%%%
We define \emph{simple parking functions} as the $\seqPF$-parking
function $(Q_j)$ satisfying $\card{Q_j} \ieq 1$ for any $j$. 
\begin{remarq}
    \label{eq:pf_simple_id_permutations}
    The simple $m$-parking functions are permutations.
\end{remarq}
\begin{ex}
    In $\PF_2(m^2)$ all parking functions are simple except
    $\boldsymbol{(}12\mid \cdot\mid\cdot\mid\cdot\boldsymbol{)}$ (see Example
    \ref{ex:pf_structures_012_3}).
\end{ex}
The interpretation of Frobenius characteristic $\ncch(\PF(\seqPF))$
(\ref{eq:ncch_PF}) reveals an easy way to extract an enumeration formula. 
The term $\S^{1^n}$ in $\ncch(\PF_n(\seqPF))$ encodes $\seqPF$-parking functions
such that $\card{Q_i} \ieq 1$ (\textit{cf}.
Remark \ref{rq:interpretation_complete_ncsf} and see
\cite[\pvexplicationSchi]{priezvirmaux}).
\begin{lemme}
    \label{lem:simpl_pf_enum}
     The \emph{simple $\seqPF$-parking functions} are enumerated by:
    \begin{align*}
    \spf(\seqPF; n) = n!\sum_{\tau \models n} \prod_{i=1}^{\ell(\tau)}
    \binom{\Xi(\seqPF; \tau, i)}{\tau_i}\,&& \text{with} \qquad \Xi(\seqPF;
    \tau, i) = \begin{dcases*}
        \seqPFn{1} & if $i = 1$\\
        \seqPFn{1 + \tau(i)} - \seqPFn{1 + \tau(i-1)} & in otherwise.
    \end{dcases*}
    \end{align*}
\end{lemme}
\begin{remarq}
    \label{rq:simple_parking_functions_trivial_automorphism}
    Simple $\seqPF$-parking functions do not have non-trivial automorphisms.
\end{remarq}
%
% \begin{remarq}
%     The formula is strongly inefficient ($2^n$ terms). The recursive version
%     is hugely more efficient:
%     \begin{align*}
%     \spf(\seqPF; n) = \sum_{k=1}^{n}
%     \binom{n}{k}\frac{\seqPFn{1}!}{(\seqPFn{1}- k)!}\spf(\psi_k; n-k)
%     \end{align*}
%     with $\psi_k(m) := \seqPFn{m+k} - \seqPFn{1}$ and $\spf(\seqPF;0)=1$ (see
%     \cite[\pvexplicationSchi]{priezvirmaux}). The solution of this recurrence is
%     Lemma \ref{lem:simpl_pf_enum}.
% \end{remarq}
%%%%%%%%%%%%%%%%%%%%%%%%%%%%
\subsubsection{Coreachability and parking functions}
\label{ssec:coreach_pf}
%%%%%%%%%%%%%%%%%%%%%%%%%%%%
From $\bij$ definition, the sets $Q_1$ and $Q_2$ (of a parking
function) encode the states which have transitions to the absorbing states.
States of $Q_1$ are non-accepting. The idea is to fix a constraint such that
it is forbidden to have states in $Q_1$.
More generally from the subsection
\ref{ssec:extended_corea_simpl_non_initial_ADFA_constraints}, we want to
consider $t+1$ constraints. The construction $\bij$ can be easily generalized to
\emph{extended ADFA}, then we remark that all transitions to all absorbing
states (extra ones and $\emptyset$) are associated to the \emph{first division
factors}.

The present parking function formalism allows us to easily consider structures
with those constraints. The first division factor depends only of $\seqPFn{1}$
and adding a (negative) constant to $\seqPF$ influence only the size of that
first factors.
\begin{thm}
    \label{thm:extended_corea_simple_niADFA_constr_pf}
    There is a bijection between \emph{extended coreachable simple
    non-initial ADFA with $t+1$ constraints} and simple $(2(m+t)^k - t -
    1)$-parking functions.
\end{thm}

%% file: sec_main.tex
In section \ref{sec:adfa} we showed that \emph{simple non-initial ADFA} $\aut$
can be described by couples $(\aut^{(i)}, \widebar{\aut}^{(i)})$ of
\emph{MADFA}, and \emph{extended coreachable simple non-initial ADFA with
constraints} (Lemma \ref{lem:niscADFA_simeq_MADFA_esADFA}). Thanks to Lemma
\ref{lem:simpl_pf_enum}, these sets of extended coreachable simple non-initial
ADFA with $t+1$ constraints have same cardinality. Finally thanks to Theorem
\ref{thm:extended_corea_simple_niADFA_constr_pf} we can enumerate the
\emph{(extended) coreachable simple ADFA}.
\begin{thm}
    The MADFA over an alphabet of $k$ symbols with $n$ states (with $i$ the
    initial state fixed) are enumerated by $\madfa(k; n)$ satisfying following
    relation:
    \begin{align*}
    \spf(2m^k-1; n) &= \sum_{t=1}^{n} \binom{n-1}{t-1} \spf(2(m+t)^k-t-1 ;n-t)
    \madfa(k; t)
    \end{align*}
    with $\madfa(k; 1) = 1$.
\end{thm}  
\begin{table}[ht]
\centering
\subfloat[][The number of (non-labelled) coreachable simple non-initial ADFA
over an alphabet of $k$ symbols: $\spf(2m^k - 1; n)/n!$.]{ 
$\scalebox{.9}{$\begin{array}{r||r|r|r|r}
n\backslash k & 1 & 2 & 3 & 4\\ \hline
0 & 1& 1& 1& 1\\
1 & 1& 1& 1& 1\\
2 & 2& 6& 14& 30\\
3 & 5& 75& 623& 4335\\
4 & 14& 1490& 59766& 1829410\\
5 & 42& 41415& 10182221& 1739056185%\\
%6 & 132& 1495806& 2740931060& 3173912127246
\end{array}$}$}
\qquad
\subfloat[][The number of (non-labelled) (quasi-)simple non-initial ADFA
over an alphabet of $k$ symbols: $\spf(2m^k; n)/n!$]{
$\scalebox{.95}{$\begin{array}{r||r|r|r|r}
n\backslash k & 1 & 2 & 3 & 4\\ \hline
0 & 1& 1& 1& 1\\
1 & 2& 2& 2& 2\\
2 & 5& 13& 29& 61\\
3 & 14& 166& 1298& 8830\\
4 & 42& 3324& 124706& 3727540\\
5 & 132& 92718& 21256346& 3543721650%\\
%6 & 429& 3354712& 5723038913& 6467723423732
\end{array}$}$}\qquad
\subfloat[][The number of (non-labelled) minimal ADFA
over an alphabet of $k$ symbols: $\madfa(k; n)/(n-1)!$]{\label{tabc}
$\scalebox{.95}{$\begin{array}{r||r|r|r|r}
n\backslash k & 1 & 2 & 3 & 4\\\hline
1 & 1& 1& 1& 1\\
2 & 2& 6& 14& 30\\
3 & 4& 60& 532& 3900\\
4 & 8& 900& 42644& 1460700\\
5 & 16& 18480& 6011320& 1220162880\\
6 & 32& 487560& 1330452032& 1943245777800\\
7 & 64& 15824880& 428484011200& 5307146859111120\\
8 & 128& 612504240& 190167920278448& 23025057433925970000\\
9 & 256& 27619664640& 111649548558856000& 149780070423407303443200\\
10 &512& 1425084870240& 84001095774695390816& 1396395902225576206029949920\\
11 &1024& 82937356685760& 78954926089415009686528&
17993790111404399137868446737600 \end{array}$}$}
\caption{Some enumeration of (extended) simple ADFA. The first values of
(\ref{tabc}) were presented in \cite{almeida2007enumeration}.}
\end{table}

%% file: doc.bbl
\begin{thebibliography}{RMA05}

\bibitem[AMR07]{almeida2007enumeration}
Almeida, Moreira, and Reis.
\newblock Enumeration and generation with a string automata representation.
\newblock {\em THEOR COMPUT SCI}, 387(2):93--102, 2007.

\bibitem[AMR08]{almeida2008exact}
Almeida, Moreira, and Reis.
\newblock Exact generation of minimal acyclic deterministic finite automata.
\newblock {\em INT J FOUND COMPUT S}, 19(04):751--765, 2008.

\bibitem[BLL98]{bergeron1998combinatorial}
Bergeron, Labelle, and Leroux.
\newblock Combinatorial species and tree-like structures.
\newblock 1998.

\bibitem[CH04]{campeanu2004maximum}
C{\^a}mpeanu and Ho.
\newblock The maximum state complexity for finite languages.
\newblock {\em J. of Automata, Languages and Combinatorics}, 9(2-3):189--202,
  2004.

\bibitem[DKS02]{domaratzki2002number}
Domaratzki, Kisman, and Shallit.
\newblock On the number of distinct languages accepted by finite automata with
  n states.
\newblock {\em J. of Automata, Languages and Combinatorics}, 7(4):469--486,
  2002.

\bibitem[Dom03]{domaratzki2003improved}
Domaratzki.
\newblock Improved bounds on the number of automata accepting finite languages.
\newblock In {\em LECT NOTES COMPUT SC}, pages 209--219. Springer, 2003.

\bibitem[Dom04]{domaratzki2004combinatorial}
Domaratzki.
\newblock Combinatorial interpretations of a generalization of the genocchi
  numbers.
\newblock {\em J. of Integer Sequences}, 7(04.3):6, 2004.

\bibitem[FS09]{flajolet2009analytic}
Flajolet and Sedgewick.
\newblock {\em Analytic combinatorics}.
\newblock cambridge University press, 2009.

\bibitem[Hop79]{hopcroft1979introduction}
Hopcroft.
\newblock {\em Introduction to automata theory, languages, and computation}.
\newblock Pearson Education India, 1979.

\bibitem[KW66]{konheim1966occupancy}
Konheim and Weiss.
\newblock An occupancy discipline and applications.
\newblock {\em SIAM J APPL MATH}, 14(6):1266--1274, 1966.

\bibitem[KY03]{kung2003goncarov}
Kung and Yan.
\newblock Goncarov polynomials and parking functions.
\newblock {\em J COMB THEORY A}, 102(1):16--37, 2003.

\bibitem[Lis06]{liskovets2006exact}
Liskovets.
\newblock Exact enumeration of acyclic deterministic automata.
\newblock {\em DISCRETE APPL MATH}, 154(3):537--551, 2006.

\bibitem[NT08]{novelli2008noncommutative}
Novelli and Thibon.
\newblock Noncommutative symmetric functions and lagrange inversion.
\newblock {\em ADV APPL MATH}, 40(1):8--35, 2008.

\bibitem[PV]{priezvirmaux}
Priez and Virmaux.
\newblock Non-commutative frobenius characteristic of generalized parking
  functions.
\newblock {\em arXiv}.

\bibitem[Rev92]{revuz1992minimisation}
Revuz.
\newblock Minimisation of acyclic deterministic automata in linear time.
\newblock {\em THEOR COMPUT SCI}, 92(1):181--189, 1992.

\bibitem[RMA05]{reis2005representation}
Reis, Moreira, and Almeida.
\newblock On the representation of finite automata.
\newblock In {\em Proc. of DCFS’05}. Citeseer, 2005.

\bibitem[SP02]{stanley2002polytope}
Stanley and Pitman.
\newblock A polytope related to empirical distributions, plane trees, parking
  functions, and the associahedron.
\newblock {\em DISCRETE COMPUT GEOM}, 27(4):603--602, 2002.

\end{thebibliography}
